\documentclass[reqno,11pt,a4paper,article,12pt]{amsart}
\usepackage{amsfonts,amssymb,fullpage}
\usepackage{xypic}
\usepackage[all]{xy}

\theoremstyle{definition}

\theoremstyle{remark}

\numberwithin{equation}{section}

\def\id{\mathrm{\mathop{id}}}
\def\opp{\mathrm{\mathop{opp}}}
\def\Ker{\mathrm{\mathop{Ker}}}
\def\Im{\mathrm{\mathop{Im}}}
\def\Tor{\mathrm{\mathop{Tor}}}

\begin{document}

\title{HOMOLOGICAL FINITENESS CONDITIONS FOR GROUPS, MONOIDS AND ALGEBRAS}
\author{Stephen J. Pride}
\maketitle


\begin{center}
\noindent{\bf\large Abstract}\\
\end{center}

Recently Alonso and Hermiller [2] introduced a homological
finiteness condition\break $bi\mbox{-}FP_n$ (here called {\it
weak} $bi\mbox{-}FP_n$) for monoid rings, and Kobayashi and Otto
[10] introduced a different property, also called $bi\mbox{-}FP_n$
(we adhere to their terminology). From these and other papers we
know that: $bi\mbox{-}FP_n \Rightarrow$ left and right $FP_n
\Rightarrow$ weak $bi\mbox{-}FP_n$; the first implication is not
reversible in general; the second implication is reversible for
group rings.  We show that the second implication is reversible in
general, even for arbitrary associative algebras (Theorem $1'$),
and we show that the first implication {\it is} reversible for
group rings (Theorem 2).  We also show that the all four
properties are equivalent for connected graded algebras (Theorem
4). A result on retractions (Theorem $3'$) is proved, and some
questions are raised.

\renewcommand{\thefootnote}{}
\footnote{Mathematics Subject Classification (MSC): 20M05, 16E40, 20J05}

\noindent{\bf\large 1. Introduction}\\

Throughout the paper $K$ will denote a fixed but arbitrary commutative ring, and
$R, S$ will denote (not necessarily commutative) rings.  All rings will have an identity, and ring
homomorphisms will be assumed to preserve the identity.\\

{\bf 1.1 Groups and monoids}\\

Let $B$ be a monoid, and let $KB$ be the corresponding monoid ring over $K$.  We have the standard
augmentation

$$
\varepsilon: KB \rightarrow K \quad b \mapsto 1 \quad (b \in B),
$$

\noindent and we can thus regard $K$ as a left $KB$-module $_BK$ with the $KB$-action via
$\varepsilon$:

$$
a.k = \varepsilon(a)k \quad (a \in KB, k \in K).
$$

\noindent Then $B$ is said to be of type {\it left}-$FP_n$ (over $K$) if there is a partial free resolution

$$
0 \leftarrow {_BK} \leftarrow P_0 \leftarrow P_1 \leftarrow \ldots \leftarrow P_n \eqno{(1)}
$$

\noindent where $P_0, P_1, \ldots, P_n$ are {\it finitely generated} free left $KB$-modules.  Similarly, we can
regard $K$ as a right $KB$-module $K_B$ via $\varepsilon$, and analogously define monoids of type
{\it right-$FP_n$} by requiring a partial resolution

$$
0 \leftarrow K_B \leftarrow P_0' \leftarrow P_1' \leftarrow \ldots \leftarrow P_n' \eqno{(2)}
$$

\noindent by finitely generated free right $KB$-modules.  These
two properties are equivalent if there is an involution $*$ on $B$
(that is a mapping $^*: B \rightarrow B$ satisfying $(bc)^* =
c^*b^*, b^{**} = b$ for all $b, c \in B$).  In particular, they
are equivalent for groups, and more generally inverse monoids (and
so in these cases we usually just use the term $FP_n$). However,
in general the left and right properties are different.  In [6],
an example is given of a monoid which is left-$FP_\infty$ (i.e.
$FP_n$ for all $n$) over $\mathbb{Z}$, but not even right-$FP_1$
over $\mathbb{Z}$ (and vice versa).

We remark that there are examples of groups which are of type
$FP_n$ over $\mathbb{Z}$ but not of type $FP_{n+1}$ over
$\mathbb{Z}$ for all $n$ [4].

We can also regard $K$ as a $(KB, KB)$-bimodule $_BK_B$ with 2-sided action
$$
a.k.a' = \varepsilon(a)k \varepsilon(a') \quad
(a, a' \in KB, \; k \in K).
$$

\noindent We can then define a finiteness condition by requiring that there exists a partial free
bi-resolution

$$
0 \leftarrow \;{_BK_B} \leftarrow F_0 \leftarrow F_1 \leftarrow \ldots \leftarrow F_n \eqno{(3)}
$$

\noindent where $F_0, F_1, \ldots, F_n$ are finitely generated
free $(KB, KB)$-bimodules. This property was introduced in [2],
and was called there {\it bi}-$FP_n$.  However, in this paper we
will call it {\it weak} bi-$FP_n$, to distinguish it from another
property discussed shortly.

As shown in [2],

$$
\mbox{left-}FP_n + \mbox{right-}FP_n \Longrightarrow \mbox{weak bi-}FP_n \eqno{(4)}
$$

\noindent(at least in the case when $K$ is a $PID$).  For, by the
K\"unneth Theorem, the tensor product over $K$ of the free partial
resolutions in (1), (2) gives a free partial bi-resolution of $_BK
\otimes_K K_B \cong {_BK_B}$.  (See [2, p.344] for details.)

There is also another ``natural" $(KB, KB)$-bimodule associated
with $KB$, namely $KB$ itself, regarded as a bimodule by left and
right multiplication.  In [10], the authors defined $B$ to be of
type {\it bi-}$FP_n$ (and we will adhere to their terminology in
this paper) if there exists a partial bi-resolution

$$
0 \leftarrow KB \leftarrow F_0 \leftarrow F_1 \leftarrow \ldots \leftarrow F_n \eqno{(5)}
$$

\noindent where $F_0, F_1, \ldots, F_n$ are finitely generated free $(KB, KB)$-bimodules.

As shown in [10],

$$
\mbox{bi-}FP_n \Longrightarrow \mbox{left-}FP_n + \mbox{right-}FP_n. \eqno{(6)}
$$

\noindent For if we apply $- \otimes_{KB}K$ to (5) then the
sequence remains exact and gives a partial resolution of $KB
\otimes_{KB}K \cong$ $_BK$ by finitely generated free left
$KB$-modules.  (See [10, p.338] for details.)

It was proved in [2] that the implication (4) is reversible for
{\it groups}.  However, in private correspondence with the second
author of [2], it emerged that no example was known to show in
general that the implication (4) is not reversible.  We will prove
that
no such example can exist.\\

\noindent{\bf Theorem 1}  {\it If a monoid is weak bi-$FP_n$ then it is both
left- and right-$FP_n$}.\\

As regards the reverse of the implication (6), an example is given
in [9] of a monoid which is left- and right-$FP_\infty$ but {\it
is not} bi-$FP_3$.  However, it has been an open question
whether (6) is reversible for {\it groups}.  We will show that this is the case.\\

\noindent{\bf Theorem 2} {\it If a group is $FP_n$ then it is bi-$FP_n$}.\\

Thus for groups the four properties weak bi-$FP_n$, left-$FP_n$, right-$FP_n$,
bi-$FP_n$ all coincide.\\

\noindent{\bf Question.} Is Theorem 2 true for inverse monoids?\\

A monoid $C$ is called a {\it retract} of a monoid $B$ if there
are monoid homomorphisms

\[
\xymatrix{ B \ar@<0.4ex>[r]^\psi & \ar@<0.4ex>[l]^\phi C } \quad
\psi\phi=\mathrm{id}_C.
\]

\noindent{\bf Theorem 3} {\it Each of the properties left-$FP_n$,
right-$FP_n$, bi-$FP_n$, weak bi-$FP_n$
is closed under retractions}.\\

(In the case of {\it groups}, this is proved in [1] for the more
general concept of {\it quasi}-retracts.)
\\

{\bf 1.2 Algebras}\\

Let $A$ be a $K$-algebra with an augmentation, that is, a $K$-algebra epimorphism

$$
\varepsilon: A \rightarrow K.
$$

\noindent Then we can regard $K$ as a left $A$-module $_AK$, or a right $A$-module $K_A$, or an
$(A, A)$-bimodule $_AK_A$ via $\varepsilon$.  Also we can regard $A$ as an $(A, A)$-bimodule
$_AA_A$ by left and right multiplication.  Then we can define $A$ to be {\it left}-$FP_n$,
{\it right}-$FP_n$, {\it weak bi}-$FP_n$, {\it bi}-$FP_n$ if there is a partial free resolution analogous
to (1), (2), (3), (5) respectively.

The argument in [2] shows that the implication (4) holds provided
$K$ is a PID and $A$ is free (or, more generally, flat) as a
$K$-module.  Also, the argument in [10] shows that the implication
(6) holds provided $A$ is free (or, more generally, projective) as
a $K$-module.
We will show that the {\it reverse} of (4) holds without restrictions.\\

\noindent{\bf Theorem 1$'$} {\it If $A$ is weak bi-$FP_n$, then
$A$ is left-$FP_n$ and
right-$FP_n$}.\\

We remark that Anick [3] showed that if $A$ can be presented as a
quotient of a finitely generated free $K$-algebra by an ideal
generated by a finite Gr\"obner base, then $A$ is left- and
right-$FP_\infty$. Recently Kobayashi [8] has improved this to
show that such an algebra is bi-$FP_\infty$.\\

A $K$-algebra $D$ is a {\it retract} of $A$ if there are
$K$-algebra homomorphisms

\[
\xymatrix{ A \ar@<0.4ex>[r]^\kappa & \ar@<0.4ex>[l]^\theta D }
\quad \kappa\theta=\mathrm{id}_D.
\]
\ Moreover, if $D$ has augmentation

$$
\varepsilon_D: D \rightarrow K
$$
\noindent then $D$ is an {\it augmented retract} if there exist
$\kappa,\theta$ as above such that $\varepsilon_{D}
\kappa=\varepsilon$ (and thus
$\varepsilon\theta=\varepsilon_{D}$).\\

\noindent{\bf Theorem 3$'$} (i) {\it The property bi-$FP_n$ for
algebras is closed under retractions.}

 (ii) {\it The properties left-$FP_n$, right-$FP_n$, weak
 bi-$FP_n$ are closed under augmented retractions.}\\

 Theorem 3 follows from this, because if a monoid $C$ is a retract
 of a monoid $B$, then the monoid algebra $KC$ is an augmented
 retract of $KB$.\\

We will also consider {\it connected graded algebras} (definitions will be given in \S5).\\

\noindent{\bf Theorem 4} {\it If a connected graded algebra is
left-$FP_n$ or
right-$FP_n$ then it is bi-$FP_n$.}\\

Thus for connected graded algebras, the four properties weak bi-$FP_n$,
left-$FP_n$, right-$FP_n$, bi-$FP_n$ all coincide.  We remark that the answer
to the following seems to be unknown:\\

\noindent{\bf Question.}  For any $n$, is there a graded algebra
of type $FP_n$ but not of type $FP_{n+1}$?\\

To prove Theorems $1'$ and $3'$ we first obtain a result
concerning what we call {\it retractive pairs} (see \S2). Theorem
$3'$ follows directly from this, and Theorem $1'$ follows by
working with the enveloping algebra $E = A \otimes_K A^{\opp}$ of
$A$ (see \S3).

The proof of Theorem 2 is given in \S4, and the proof of Theorem 4 is given in \S5.\\

\noindent{\bf Acknowledgements.}  I thank Peter Kropholler and Alexandro Olivares for
helpful discussions.\\

{\bf\large 2. Retractive pairs and the Property $FP_n$}\\



The following consequence of the Generalised Schanuel Lemma is useful [5, p 193].\\

\noindent{\bf Lemma 1 } {\em Let $M$ be a left $R$-module.  For $n \geq 0$, if
$$
0 \stackrel{\partial_{-1}}{\longleftarrow} M \stackrel{\partial_0}{\longleftarrow}
P_0 \longleftarrow P_1 \ldots
\stackrel{\partial_{n-1}}{\longleftarrow} P_{n-1}
$$
is a partial free resolution of $M$ of length $n-1$, with $P_0,
\ldots, P_{n-1}$ finitely generated free modules, then $M$ is of
type $FP_n$ if and only if
$\ker \partial_{n-1}$ is finitely generated}.\\

Suppose we have ring homomorphisms


\[
\xymatrix{
R \ar@<0.4ex>[r]^\rho &
\ar@<0.4ex>[l]^\iota S }
\quad \rho\iota=\mathrm{id}_S.
\]

\noindent {\it ie}, $\rho$ is a {\em retraction} of $R$ onto $S$,
with {\em section} $\iota$.

A (left) {\em retractive pair} consists of a left $R$-module $M$,
a left $S$-module $L$, and abelian group homomorphisms

$$
\xymatrix{
M \ar@<0.4ex>[r]^{\alpha^+} &
\ar@<0.4ex>[l]^{\alpha^-} } L
\quad \alpha^+ \alpha^- = \id_L,
$$

\noindent where $\alpha^+$ is an $R$-module homomorphism (regarding $L$ as an $R$-module via $\rho$),
and $\alpha^-$ is an $S$-module homomorphism (regarding $M$ as an
$S$-module via $\iota$).  We will denote such a retractive pair by $M \alpha L$.

Retractive pairs form a category, where a {\em mapping}
$$
(\phi, \psi): M \alpha L \longrightarrow M' \beta L'
$$
consists of an $R$-module homomorphism
$$
\phi: M \longrightarrow M'
$$
\noindent and an $S$-module homomorphism
$$
\psi: L \longrightarrow L',
$$
such that the diagram

$$
\xymatrix{
M\ar[rr]^{\phi}\ar@<-0.5ex>[d]_{\alpha^+} && M'\ar@<-0.5ex>[d]_{\beta^+} \\
L\ar[rr]^{\psi}\ar@<-0.5ex>[u]_{\alpha^-} && L'\ar@<-0.5ex>[u]_{\beta^-}
}
$$

\noindent commutes.  It is then easily checked that
$\alpha^+(\Ker \phi) \subseteq \Ker \psi$ and
$\alpha^-(\Ker \psi) \subseteq \Ker \phi$, so by restriction, we get the
retractive pair $\Ker(\phi, \psi)$:

$$
\xymatrix{
\Ker \phi
\ar@<0.4ex>[r]^-{\alpha^+} &
\ar@<0.4ex>[l]^-{\alpha^-} }
\Ker \psi.
$$

\noindent Similarly, we get the retractive pair $\Im(\phi, \psi)$:

$$
\xymatrix{
\Im \phi
\ar@<0.4ex>[r]^-{\beta^+} &
\ar@<0.4ex>[l]^-{\beta^-} }
\Im \psi.
$$

\noindent{\bf Proposition 1}  {\it Let $M \alpha L$ be a
retractive pair.  If $M$ is of type $FP_n$ then there is a
sequence}
$$
0 \longleftarrow M \alpha L
\stackrel{(\partial_0, \delta_0)}{\longleftarrow} P_0 \beta_0 F_0
\stackrel{(\partial_1, \delta_1)}{\longleftarrow} P_1 \beta_1
F_1 \longleftarrow \ldots
\stackrel{(\partial_n, \delta_n)}{\longleftarrow} P_n \beta_n F_n
$$
{\it with $P_i, F_i$ finitely generated free $R$-modules, $S$-modules respectively
$(0 \leq i \leq n)$, $\Im(\partial_0, \delta_0) = M \alpha L$, and
$\Im(\partial_{i+1}, \delta_{i+1}) = \Ker(\partial_i, \delta_i)\;
(0 \leq i < n)$.

In particular,
$$
0 \longleftarrow L
\stackrel{\delta_0}{\longleftarrow} F_0
\stackrel{\delta_1}{\longleftarrow} F_1
\longleftarrow \ldots
\stackrel{\delta_n}{\longleftarrow} F_n
$$
is a partial free resolution of $L$, so $L$ is of type} $FP_n$.\\

\noindent{\bf Proof.} Suppose $M$ is of type $FP_0$ (i.e. finitely generated), and let
$\{m_e: e \in \mathbf{e}\}$ be a finite set of $R$-module generators for $M$.  Then
$\{\alpha^+(m_e): e \in \mathbf{e}\}$ is a set of $S$-module generators for $L$.
Let $P_0$ be the free $R$-module (of rank $2|\mathbf{e}|$)

$$
P_0 = \left(\displaystyle{\oplus_{e \in \mathbf{e}}}Re \right) \oplus
\left(\displaystyle{\oplus_{e \in \mathbf{e}}} Re' \right),
$$
\smallskip

\noindent and let $F_0$ be the free $S$-module (of rank $|\mathbf{e}|$)

$$
F_0 = \oplus_{e \in \mathbf{e}} S \overline{e}.
$$
Then
$$
P_0
\xymatrix{
\ar@<0.4ex>[r]^{\beta^+_0} &
\ar@<0.4ex>[l]^{\beta^-_0} F_0}
\qquad
\beta^+_0(e) = \overline{e},
\beta^+_0(e') = 0,
\beta^-_0(\overline{e}) = e \;
(e \in \mathbf{e})
$$

\noindent is a retractive pair.  We have the surjective $R$-module homomorphism
$$
\partial_0: P_0 \rightarrow M, \quad
\partial_0(e) = \alpha^- \alpha^+(m_e),
\partial_0(e') = m_e - \alpha^- \alpha^+(m_e) \; (e \in \mathbf{e})
$$

\noindent and the surjective $S$-module homomorphism
$$
\delta_0: F_0 \rightarrow L, \quad \delta_0(\overline{e}) = \alpha^+(m_e) \; (e \in \mathbf{e}),
$$

\noindent and it is easily checked that $(\partial_0, \delta_0)$ is a mapping of retractive pairs.

Let $M_1 \beta_0 L_1 = \Ker(\partial_0, \delta_0)$.  By Lemma 1, if $M$ is if type
$FP_1$, then $M_1$ is finitely generated.  We can then repeat the above procedure
to obtain a finitely generated free retractive pair $P_1 \beta_1 F_1$ and a surjective map
$$
P_1 \beta_1 F_1 \rightarrow M_1 \beta_0 L_1.
$$
Composing this with the inclusion of $M_1 \beta_0 L_1$ into
$P_0 \beta_0 F_0$ we obtain a mapping
$$
(\partial_1, \delta_1): P_1 \beta_1 F_1 \rightarrow P_0 \beta_0 F_0.
$$
Continuing in this way, after $n+1$ steps we get the required sequence.\\

{\bf\large 3. Bi-Resolutions and Enveloping Algebras}\\

Recall that an $(A, A)$-bimodule $M$ is an abelian group on which $A$ acts on the left and right, with
the condition that: $(am)b=a(mb)$ for all
$a, b \in A, m \in M$; $km=mk$ for all $m \in M, k \in K$.

The $(A, A)$-bimodule $A \otimes_K A$ (with bi-action given by
$a.(u \otimes v).b = au \otimes vb$ for all $a, b, u, v \in A$) is free on the
generator $1 \otimes 1$.  Thus a direct sum of $r$ copies of $A \otimes_K A$ is a
{\em free $(A, A)$-bimodule of rank} $r$ ($r$ may be infinite).

An $(A, A)$-bimodule $M$ is said to be of {\em type bi}-$FP_n$ if there is a partial resolution
$$
0 \leftarrow M \leftarrow F_0 \leftarrow F_1 \leftarrow \ldots \leftarrow F_n \eqno{(7)}
$$
where $F_0, F_1, \ldots, F_n$ are {\em finitely generated} free
$(A, A)$-bimodules.

It is obvious that Proposition 1 extends to bimodules.\\

\noindent{\bf Proof of Theorem 3$'$.} Let $A$, $D$ be as in the
paragraph before the statement of Theorem 3$'$.

(i) Regarding $A$, $D$ as ($A$, $A$)-, ($D$, $D$)-bimodules,
respectively, ${\xymatrix{ A \ar@<0.3ex>[r]^\kappa &
\ar@<0.3ex>[l]^\theta D }}$ is a retractive pair.

(ii) Regarding $K$ as a left (respectively, right, bi) $A$-module,
and a left (respectively, right, bi) $D$-module, ${\xymatrix{ K
\ar@<0.3ex>[r]^\id & \ar@<0.3ex>[l]^\id K }}$ is a retractive
pair.\\

Recall that for a $K$-algebra $A$, there is the {\em opposite} algebra
$A^{opp}$.  This has the same underlying set as $A$, and the same addition and
scalar multiplication.  The product of two elements $a, b$ (in that order) in
$A^{opp}$ is defined to be the product $ba$ in $A$.  When regarding an element $a$ of
$A$ as an element of $A^{opp}$, we will denote it by $a^{opp}$.

The {\em enveloping algebra} of $A$ is the tensor product
$E = A \otimes_K A^{opp}$, with multiplication defined by
$$
(a \otimes b^{opp}) (c \otimes d^{opp}) = ac \otimes b^{opp} d^{opp} =
ac \otimes db \quad (a,b, c, d \in A). \eqno{(8)}
$$
There is the induced augmentation
$$
\varepsilon_E: E \rightarrow K \quad  \varepsilon_E(a \otimes
b^{opp}) = \varepsilon (a) \varepsilon (b)  \quad (a, b \in A).
$$

If $M$ is an $(A, A)$-bimodule we can regard it as a left $E$-module (denoted
$\mathcal{E}(M)$) with $E$-action given by

$$
(a \otimes b^{opp})m = amb \quad (a, b \in A, m \in M).
$$
Also, if $\phi: M \rightarrow M'$ is a bimodule homomorphism, then it can be regarded as a left
$E$-module homomorphism $\mathcal{E}(\phi): \mathcal{E}(M) \rightarrow \mathcal{E}(M')$.  Then
$\mathcal{E}$ is an (exact) functor.  This functor has an inverse $\mathcal{A}$, where for a left
$E$-module $N$, $\mathcal{A}(N)$ is $N$ regarded as an $(A, A)$-bimodule, with left and right
actions given by
$$
an = (a \otimes 1)n, \; nb = (1 \otimes b^{opp})n \quad
(a, b \in A, n \in N).
$$

It is easily shown that $\mathcal{E}(A \otimes_K A)$ is $E$ acting on itself by left
multiplication.  In other words, $\mathcal{E}(A \otimes_K A)$ is a free left $E$-module of rank 1.  Thus,
if $F$ is a free $(A, A)$-bimodule of rank $r$, then $\mathcal{E}(F)$ is a free left $E$-module of rank
$r$.  Applying $\mathcal{E}$ to a partial resolution as in (7), we thus see that if $M$ is
an $(A, A)$-bimodule of type bi-$FP_n$, then $\mathcal{E}(M)$ is of type $FP_n$.
By considering the inverse functor $\mathcal{A}$, the converse is also true.  Thus we have:\\

\noindent{\bf Lemma 2} {\em An $(A, A)$-bimodule $M$ is of type bi-$FP_n$ if and only if
$M$, regarded as a left $E$-module, is of type} $FP_n$.\\

\noindent{\bf Proof of Theorem $1'$.} \noindent Regarding $K$ as
an $(A, A)$-bimodule via $\varepsilon$, $\mathcal{E}(K)$ is easily
seen to be $K$ regarded as a left $E$-module via $\varepsilon_E$.
Thus, by Lemma 2, $A$ is weak bi-$FP_n$ if and only if $E$ is
left-$FP_n$. Then, since $A$ is an augmented retract of $E$ under
the maps
$$
E \xymatrix{ \ar@<0.4ex>[r]^{\rho} & \ar@<0.4ex>[l]^{\iota}
A}\quad \rho(a \otimes b^{opp}) = \varepsilon(b)a, \; \iota(a) = a
\otimes 1 \; (a, b \in A),\; \rho \iota = id_A,
$$
if $E$ is left-$FP_n$ then so is $A$, by Theorem 3$'$(ii).\\

{\bf\large 4. Proof of Theorem 2}\\

For a group $G$ we define a functor -- $\hat{\otimes}KG$ from the category of left
$KG$-modules to the category of $(KG, KG)$-bimodules as follows.
For $M$ a left $KG$-module, $M \hat{\otimes}KG$ is the tensor product
$M \otimes_K KG$ with bi-$KG$-action given by
$$
g \cdot (m \otimes x) \cdot h = gm \otimes g xh \quad
(g,h,x \in G, \; m \in M).
$$

\noindent For a left $KG$-module homomorphism
$\alpha: M_1 \rightarrow M_2$ we define
$\hat{\alpha}: M_1 \hat{\otimes} KG \rightarrow M_2 \hat{\otimes} KG$ to be
$$
\alpha \otimes \id_{KG}: M_1 \otimes_K KG \rightarrow M_2 \otimes_K KG,
$$

\noindent regarded as a $(KG, KG)$-bimodule homomorphism.\\

\noindent{\bf Lemma 3} (i) $K \hat{\otimes} KG$ {\it is isomorphic
to $KG$ (regarded as a $(KG, KG)$-bimodule by left and right multiplication)

(ii) If $P$ is a free left $KG$-module of rank $r$, then
$P \hat{\otimes} KG$ is a free $(KG, KG)$-bimodule of rank} $r$.\\

\noindent{\bf Proof.} (i) As an abelian group, $K \hat{\otimes}KG$
is just $K \otimes_K KG$, which is isomorphic to $KG$ by the
isomorphism

$$
\theta: K \otimes_K KG \rightarrow KG \qquad
k \otimes x \mapsto kx \; (k \in K, x \in G).
$$

It is easily checked that for $g, h, x \in G, \; k \in K$

$$
\begin{array}{lll}
\theta(g \cdot (k \otimes x) \cdot h)&=& g \theta (k \otimes x) h.
\end{array}
$$

(ii) Since $P$ is the direct sum of $r$ copies of $KG$, it suffices to show that
$KG \hat{\otimes} KG$ is a free $(KG, KG)$-module of rank 1.  The free $(KG, KG)$-bimodule of rank 1 is
$KG \otimes_K KG$ with action
$$
x(g \otimes h)y = xg \otimes hy \quad (x,y, g, h \in G).
$$

\noindent As a $K$-module, $KG \otimes_K KG$ is free with basis
$g \otimes h \; (g, h \in G)$.  For convenience write
$g \hat{\otimes} h$ for $g \otimes h$ when considered as an element
of $KG \hat{\otimes} KG$.


Since $KG \otimes_K KG$ is free on ${1 \otimes 1}$, we get a bi-module
homomorphism
$$
\beta: KG \otimes_K KG \rightarrow KG \hat{\otimes}KG \quad
1 \otimes 1 \mapsto 1 \hat{\otimes}1.
$$

\noindent Thus $\beta(g \otimes h) = \beta(g. (1 \otimes 1.h) =
g \hat{\otimes} gh \; (g, h \in G)$.  Also, we have a $K$-module homorphism
$$
\alpha: KG \hat{\otimes} KG \rightarrow KG \otimes_K KG \quad
g \hat{\otimes} h \mapsto g \otimes g^{-1}h \; (g, h \in G).
$$

\noindent It is easily checked that
$\alpha$ and $\beta$ are mutually inverse (as $K$-maps), so $\beta$ is a $(KG, KG)$-bimodule isomorphism.\\

To prove Theorem 2, suppose $G$ is left-$FP_n$ over $K$. Then there is a partial free resolution
as in (1).  Tensoring by $- \otimes_K KG$ is exact, since $KG$ is free as a $K$-module,
so we obtain the exact sequence
$$
0 \leftarrow K \otimes_K KG \leftarrow P_0 \otimes_K KG \leftarrow P_1
\otimes_K KG \leftarrow \cdots \leftarrow P_n \otimes_K KG.
$$

\noindent Regarding this as a sequence of $(KG, KG)$-bimodules and using Lemma 3,
we obtain an exact sequence
$$
0 \leftarrow KG \leftarrow F_0 \leftarrow F_1 \leftarrow \cdots \leftarrow F_n
$$
where $F_0, F_1, \cdots, F_n$ are finitely generated free $(KG, KG)$-bimodules,
so $G$ is bi-$FP_n$.\\

{\bf\large 5.  Connected Graded Algebras and Proof of Theorem 4}\\

Suppose $K$ is a field, and that $A$ is a {\it graded} $K$-algebra.  Thus
$A$ is a direct sum $\oplus_{i \geq 0} A_i$ of $K$-modules such that
$A_i A_j \subseteq A_{i+j}\;(0 \leq i, j)$.  Elements of $A_i$ are said to be of
{\it degree} $i$.  In the context of graded algebras, modules will also be graded.
Thus a left $A$-module is a directed sum $M = \oplus_{i \geq 0} M_i$ of $K$-modules
such that $A_i M_j \subseteq M_{i+j} \; (0 \leq i, j)$.  Right modules and bimodule
are defined analogously.  A module is concentrated in dimension $n$ if $M_i = 0$
for $i \neq n$.  We regard $K$ as a graded module concentrated in degree $0$.  A mapping
$\phi: M \rightarrow L$ of left modules consists of a family $\phi_i: M_i \rightarrow L_i \quad (i \geq 0)$ of
$K$-maps such that $\phi_{i+j}(a_i m_j) = a_i \phi_j(m_j) \; (a_i \in A_i,
m_j \in M_j, i, j \geq 0)$.

We will assume that $A$ is {\it connected}, that is, $A_0$ has basis the identity
$1_A$.  Then we have the standard augmentation
$$
\varepsilon: A \rightarrow K \quad
1_A \mapsto 1_K, \quad
A_i \rightarrow 0 \; (i > 0)
$$

\noindent with kernel $A^+ = \oplus_{i > 0} A_i$.

The opposite algebra $A^{\opp}$ is also a connected graded algebra with the
same grading, and so $E = A \otimes_K A^{\opp}$ inherits a connected graded
algebra structure with grading
$$
E_i = \oplus_{p+q=i} A_p \otimes_K A_q
$$

\noindent and multiplication as in (8).\\

\noindent{\bf Remark} If $M$ is a graded left $A$-module then there are two
possible definitions of $FP_n$, according to whether we consider free resolutions

$$
0 \longleftarrow M \stackrel{\partial_0}{\longleftarrow} P_0
\stackrel{\partial_1}{\longleftarrow} P_1 \longleftarrow \ldots
\stackrel{\partial_i}{\longleftarrow} P_i \longleftarrow \ldots \eqno{(9)}
$$

\noindent where the $P_i$'s are free modules which are graded and
the $\partial_i$'s are graded maps (``$FP_n$ in the graded
sense"), or whether we just consider ungraded resolutions
(``$FP_n$ in the ungraded sense").  Clearly, if $M$ is $FP_n$ in
the graded sense then it is $FP_n$ in the ungraded sense.  The
converse is also true.  For if $M$ is $FP_n$ in the ungraded sense
then applying $K \otimes_A-$ to an ungraded resolution (9) with
$P_0, P_1, \ldots, P_n$ finitely generated, we see that
$\Tor^A_i(K,M)$ is a finitely generated $K$-module for $0 \leq i
\leq n$.  Now associated with any graded module is a canonical
(minimal) graded resolution [7], where the $i$th term is $A
\otimes_K \Tor^A_i(K,M)$, so $M$ is of type $FP_n$ in the graded
sense.  Similar remarks
also hold for right modules, and bimodules.\\

\noindent{\bf Proposition 2} {\it Let $M$ be an $(A, A)$-bimodule
which as a right $A$-module is free.  If the left $A$-module $M
\otimes_A K$ is of type $FP_n$, then the bimodule $M$ is of type
$FP_n$.  (An analogous result holds if we interchange left
and right}.)\\

Taking $M = _AA_A$ in Proposition 2 we obtain Theorem 4.\\

{\bf Proof.} Consider the left $E$-module $\mathcal{E}(M)$.  By
standard theory [7], there is a unique (up to isomorphism) minimal
resolution

$$
0 \longleftarrow \mathcal{E}(M) \stackrel{\partial_0}{\longleftarrow} P_0
\stackrel{\partial_1}{\longleftarrow} P_1 \longleftarrow \ldots
\stackrel{\partial_i}{\longleftarrow} P_i \longleftarrow \ldots
$$

\noindent where $P_i$ is a free $E$-module, and $\Ker \partial_i \subseteq E^+ P_i$.
Applying the functor $\mathcal{A}$, we then get a bi-resolution

$$
0 \longleftarrow M \stackrel{\overline{\partial}_0}{\longleftarrow }
\overline{P}_0 \stackrel{\overline{\partial}_1}{\longleftarrow}
\overline{P}_1 \longleftarrow \ldots
\stackrel{\overline{\partial}_i}{\longleftarrow}
\overline{P}_i \longleftarrow \ldots .
$$

\noindent where $\overline{P}_i$ is a free $(A, A)$-bimodule and
$\Ker \overline{\partial_i} \subseteq \mathcal{A}(E^+P_i) =
A^+ \overline{P}_i + \overline{P}_i A^+$.  Note that since a free
$(A, A)$-bimodule is also free as a right
$A$-module, the above sequence is also a free right resolution of $M$ regarded as
a right $A$-module.

Now applying  $- \otimes_A K$ we obtain

$$
0 \longleftarrow M \otimes_A K \stackrel{\overline{\partial}_0 \otimes 1}{\longleftarrow}
\overline{P}_0 \otimes_A K \stackrel{\overline{\partial}_1 \otimes 1}{\longleftarrow}
\overline{P}_1 \otimes_A K
\longleftarrow  \ldots \stackrel{\overline{\partial}_i \otimes 1}{\longleftarrow}
\overline{P}_i \otimes_A K \longleftarrow \ldots .
$$

\noindent Then for $i \geq 1$

$$
\begin{array}{lll}
\frac{\Ker \overline{\partial}_i \otimes 1}{\Im \overline{\partial}_{i+1} \otimes 1} &=&
\Tor^A_i (M, K)\\[5pt]
&=& 0 \; \mbox{ since } M \mbox{ is free as a right } A\mbox{-module}.
\end{array}
$$

\noindent So the sequence is exact. Since $\Ker(\overline{\partial}_i \otimes 1)
\subseteq A^+ (\overline{P}_i \otimes 1)$ the sequence is, in fact, the minimal resolution
of $M \otimes_A K$.  Thus if $M \otimes_A K$ is of type $FP_n$, then the left $A$-modules
$\overline{P}_j \otimes_A K \; (0 \leq j \leq n)$ are finitely generated.  Thus the
$(A, A)$-bimodules
$$
(\overline{P}_j \otimes_A K) \otimes_K A \cong \overline{P}_j \quad
(0 \leq j \leq n)
$$

\noindent are finitely generated, so the bimodule $M$ is of type $FP_n$.

\vspace{0.10cm}

\begin{tabular}{lll}
Address of author: & \\
Department of Mathematics&\\
University of Glasgow &\\
15 University Gardens &\\
Hillhead &\\
Glasgow &\\
G12 8QW &\\
&\\
e-mail address: & \\
sjp@maths.gla.ac.uk &\\
\end{tabular}

\end{document}